\documentclass[11pt]{amsart} 
\usepackage{graphicx}
\usepackage{vmargin}
\usepackage{subfig}

\def\perm{\pi}
\def\trans{\tau}

\def\indset{\mathcal{I}}

\def\ext#1{#1^*} 
\def\rest#1{#1_*} 
\def\up#1{#1^\uparrow} 
\def\down#1{#1_\downarrow} 

\DeclareMathOperator{\sym}{S}
\DeclareMathOperator{\cyc}{C}
\DeclareMathOperator{\rand}{\texttt{rand}}
\DeclareMathOperator{\swap}{\texttt{swap}}

\DeclareMathOperator{\iso}{iso}
\DeclareMathOperator{\inc}{\iota}

\newcommand{\algfont}{\texttt}

\newcommand{\Em}[1]{\textbf{\textit{#1}}}

\newcommand{\Algorithm}[4]{ 
\begin{tabbing}
xxxx\=xxxx\=xxxx\=xxxx\=xxxx\=xxxx\=xxxx\= \kill
\textbf{algorithm } \algfont{#1} \\
\> {\textit{\textbf{Input: }} #2}\\
\> {\textit{\textbf{Output: }} #3} \\
\textbf{begin} \\
#4
\textbf{end}
\end{tabbing}
}

\numberwithin{equation}{section}
\newtheorem{theorem}{Theorem}[section]
\newtheorem{proposition}[theorem]{Proposition}

\theoremstyle{definition}
\newtheorem{definition}{Definition}[section]

\title{Random and exhaustive generation of permutations and cycles}
\author{Mark C. Wilson}
\address{Department of Computer Science, University of Auckland,
Private Bag 92019 Auckland, New Zealand}
\email{mcw@cs.auckland.ac.nz}

\thanks{Thanks to Hosam Mahmoud and Philippe Flajolet for useful discussions.}

\begin{document}

\begin{abstract}
In 1986 S.~Sattolo introduced a simple algorithm for uniform random
generation of cyclic permutations on a fixed number of symbols. This
algorithm is very similar to the standard method for generating a random
permutation, but is less well known.

We consider both methods in a unified way, and discuss
their relation with exhaustive generation methods. We analyse several
random variables associated with the algorithms and find their grand
probability generating functions, which gives easy access to moments and 
limit laws.

\end{abstract} 

\maketitle

\section{The algorithms}
\label{sec:algo}

\subsection*{Basic notation}
\label{ss:notation}

For each $n\geq 1$, we denote by  $\sym_n$ the \Em{symmetric group} (set of all
permutations under the operation of composition) on the set $[n]:=\{1,
\dots, n\}$ . Suppose that $n\geq 2$. For each  $\perm \in \sym_{n - 1}$,
let $\ext{\perm}$ be its \Em{extension} to $\sym_n$, by definition the
element of $\sym_n$ that fixes $n$ and agrees with $\perm$ on $[n-1]$. The
map $\ext{}$ is injective on each $\sym_n$. Dually, for an element $\rho
\in \sym_n$ that fixes $n$, $\rest{\rho}$ is the \Em{restriction} of $\rho$
to  $\sym_{n - 1}$. The map $\rest{}$ is onto each $\sym_{n - 1}$ and
$\ext{}$ followed by $\rest{}$ is the identity on $\sym_{n - 1}$. Thus we
may consider that $\sym_n \subset \sym_{n+1}$ and let $\sym$ be the union
of all $\sym_n$ (formally, we consider the direct limit induced by the
natural inclusion maps $[n] \to [n+1]$). Each element $\perm$ of $\sym$ 
belongs to a maximal $\sym_n$, where $n$ is the largest integer moved by 
$\perm$; we define $n(\perm)$ to be this value of $n$.

The \Em{action} of $\perm\in \sym$ on $i$ is denoted by $i \perm$. We use the
standard representation as words throughout; the element $\perm$ of $\sym_n$ is
written as the word $1 \perm \cdots n \perm$. Let $\cyc_n$ be the set of 
$n$-\Em{cycles} of $\sym_n$ (recall that an element of $\sym$ is a $k$-cycle 
if and only
if its action has a single nontrivial orbit, and this orbit has size $k$). When
$n = 1$, our convention is that $\cyc_n = \sym_n$. A $2$-cycle is called a
\Em{transposition}, and we denote by $\trans(i, j)$ the transposition that
exchanges $i$ and $j$ and fixes all other symbols.

Finally, for $\perm\in \sym_n$, we define $q(\perm) = n  \perm^{-1}$, and 
$Q(\perm) = n \perm = q (\perm^{-1})$.

\subsection*{Random generation}
\label{ss:algo}

Pseudocode for the two algorithms discussed below is shown in 
Figure~\ref{fig:pseudo}. Note that for notational convenience we will consider permutations 
on $\{0, \dots, n-1\}$ instead of $\{1, \dots, n\}$ in that figure, as well as in
 section~\ref{ss:exhaustive}.

The standard algorithm \cite[3.4.2, Algorithm P]{knuth;taocp2} for
uniformly generating a random permutation of $[n]$  is as follows. Start
with the identity permutation. There are $n - 1$ steps. At the $i$th step,
a random position $j$ is chosen uniformly from $[n - i + 1]$ and the
current element in position $j$ is swapped with the element at position $n
- i + 1$. Example: the permutation $25314 \in \sym_5$ is formed by
choosing $j = 4, 1, 3, 1$ in that order. Knuth attributes this algorithm to
R. A. Fisher and F. Yates \cite{fisher-yates}, and a computer implementation was 
given by Durstenfeld \cite{durstenfeld}; it is often called the
``Fisher-Yates shuffle" or the ``Knuth shuffle".

S. Sattolo \cite{sattolo} introduced a very similar algorithm for uniform
random generation of an element of $\cyc_n$.  The only difference in the
algorithm is that the possibility $j = n - i + 1$ is disallowed; $j$ is
chosen uniformly at random from $[n - i]$.

\begin{figure}
\begin{minipage}{0.5\linewidth}
\Algorithm{fisher-yates}{positive integer $n$}{permutation $\perm \in \sym_n$}
{
\textbf{for} $i$ \textbf{from}  $0$ \textbf{to} $n-1$ \textbf{do}\\
\> $\perm[i]: = i$\\
\textbf{if} ($n > 1$) \textbf{then} \\
\> \textbf{for} $k$ \textbf{from} $n-1$ \textbf{downto} $1$ \textbf{do} \\
\> \> $j : = \rand(0..k)$ \\
\> \> $\perm \gets \swap(\perm, r, k)$ \\
return($\perm$)\\
}
\end{minipage}
\hspace{1.0cm}
\begin{minipage}{0.5\linewidth}
\Algorithm{sattolo}{positive integer $n$}{permutation $\perm \in \cyc_n$}
{
\textbf{for} $i$ \textbf{from}  $0$ \textbf{to} $n-1$ \textbf{do}\\
\> $\perm[i]: = i$\\
\textbf{if} ($n > 1$) \textbf{then} \\
\> \textbf{for} $k$ \textbf{from} $n-1$ \textbf{downto} $1$ \textbf{do} \\
\> \> $j : = \rand(0..k-1)$ \\
\> \> $\perm \gets \swap(\perm, r, k)$ \\
return($\perm$)\\
}
\end{minipage}
\caption{The random generation algorithms}
\label{fig:pseudo}
\end{figure}

\subsection*{Algebraic description}
\label{ss:algebra}

In terms of multiplication in the symmetric group $\sym_n$, the description
of the Fisher-Yates algorithm is as follows. For each $i, j \in [n]$, let $\trans(i, j)$ denote the
transposition that exchanges $i$ and $j$. Then each algorithm starts with
the identity permutation $\perm$. At each value of $k$, the statement $\perm \gets
\perm \trans(r \perm, j \perm)$ is executed. In terms of
multiplication on the left, we have the following. Observe that the
transposition $\trans (r \perm, j  \perm)$ equals the conjugation
$\perm^{-1} \trans(r, j) \perm$. Thus we can rewrite the step above as 
$\perm \gets \trans(r, j) \perm$.

\begin{definition}
We call a product $\trans_1 \cdots \trans_k$ of transpositions a
\Em{triangular product} if for each $i\leq k$, $\trans_i$ exchanges
$i + 1$ and some $j \leq i + 1$. The product is called a \Em{strict triangular product} if 
always $j \leq i$. 
\end{definition}

It follows directly from the above discussion that each execution of the Fisher-Yates (respectively 
Sattolo's) algorithm yields a triangular (respectively strictly triangular) product of 
$n - 1$ transpositions in $\sym_n$. Furthermore these maps are 1--1. To see this, note that 
given $\perm$ we can reconstruct $\trans_i$ for each $i$. This is because $\trans_i$ 
fixes all elements greater than $i+1$, so that $n \cdot \trans_{n-1} = n \cdot \perm$. 
This determines the transposition $\trans_{n-1}$, and the result follows by induction on $n$. 

We may therefore define maps $\up{}: \perm_{n - 1} \times [n-1] \to \perm_n$
and $(\down{}, q): \perm_n \to \perm_{n-1} \times [n-1]$ by 

\begin{align*}
\up{(\trans_1\dots \trans_{n-2}, q)} & = \trans(q, n) \trans_1\dots \trans_{n-2} = \trans_1\dots \trans_{n-2} \trans(n, q \cdot \trans_{q-1})  \\
\down{(\trans_1 \dots \trans_{n-1})} & = \trans_1 \dots \trans_{n-2} 
\end{align*}
Note that $\down{(\perm^{-1})} = (\down{\perm})^{-1}$. Note also that we could also define $\up{}$ 
and $\down{}$ directly without reference to the triangular representation:

\begin{align*}
\up{(\perm, q)} & = \trans(n(\perm) + 1, q) \ext{\perm} = \ext{\perm} \trans(n(\perm) + 1, q \cdot \perm)) \\
\down{\perm} & = \rest{(\trans(n(\perm), q(\perm)) \perm)} =
\rest{(\perm \trans(n(\perm), Q(\perm)))}.
\end{align*}

\begin{proposition}
\label{prop:decomp}
For $n\geq 2$, the maps $\, \up{}: \perm_{n - 1} \times [n-1] \to \perm_n$
and $\, (\down{}, q): \perm_n \to \perm_{n-1} \times [n-1]$ defined  above are mutually inverse bijections.
Furthermore each restricts to $\cyc$ and the restrictions are also mutually inverse bijections.
\end{proposition} 
\begin{proof}
Note that if $\perm'$ denotes $\trans(n(\perm) + 1, q) \ext{\perm}$, then $q(\perm') = q$ and $n(\perm') = n(\perm) + 1$. Thus the composition of the maps in either order is the identity. 
Suppose that $\sigma \in \cyc_{n-1}$ and $1\leq q < n-1$. Then $\rho:=\up{\sigma}$ has the property that $n \cdot \rho^i = q \cdot \sigma^i$ for $1 \leq i \leq {n-1}$, and hence never equals $n$. Thus 
$\rho$ is an $n$-cycle.
\end{proof}

By iteration this yields a map $\iso_n: \sym_n \to [n]!:=[n] \times [n-1]
\times \dots \times [2] \times [1]$ taking a permutation to the sequence of
positions $j$ made in the execution of the Fisher-Yates algorithm.
 Sattolo's algorithm fits nicely into this picture. For each 
$i$, there is a natural inclusion map $[i] \to [i+1]$. The product of these gives a map
$\inc_{n-1} : [n-1]! \to  [n]!$. 

We summarize the above result in a proposition.

\begin{proposition}
The following sets are in bijection via the correspondences described above. 
\begin{enumerate}
\item the set of possible outputs of the Fisher-Yates algorithm;
\item $\sym_n$;
\item the set of triangular decompositions of length $n - 1$ in $\sym_n$;
\item the set $[n] \times [n-1] \times \dots \times [2] \times [1]$.
\end{enumerate}
Furthermore, the measure induced on $\sym_n$ by the Fisher-Yates algorithm is uniform.

The following sets are in bijection via the correspondences described above.
\begin{enumerate}
\item the set of possible outputs of Sattolo's algorithm;
\item $\cyc_n$;
\item the set of strict triangular decompositions of length $n - 1$ in $\sym_n$;
\item the set $[n-1] \times [n-2] \times \dots \times [2] \times [1]$.
\end{enumerate}
Furthermore, the measure induced on $\cyc_n$ by Sattolo's algorithm is uniform.
\end{proposition}

Note that the initial subproduct of length $i$ of a (strict) triangular
product $\trans_1 \dots \trans_{n-1}$ of transpositions in $\sym_n$ is
itself a (strict) triangular product of transpositions in $\sym_{i+1}$, and
hence an $(i+1)$-permutation/cycle. At each stage, forming the next partial
product involves inserting $i+1$ into the current permutation/cycle. 
This gives an algorithm for forming a uniformly random permutation/cycle of a random 
length; simply form such a (strict) triangular product with length chosen 
according to the desired distribution. Clearly the distribution conditioned on the 
length is uniform.

\subsection*{Exhaustive generation}
\label{ss:exhaustive}

There are obvious deterministic versions of the above algorithm. Instead of
randomly choosing the transpositions, we simply run through all such
transpositions systematically. Every method of generating all elements of
$[n]!$ (and the corresponding unranking function) can be transferred via
the encoding above to a method for generating all elements of $\sym_n$
or $\cyc_{n+1}$.

A common way of enumerating a combinatorial class is to use an
\Em{incremental method}, where each object is generated from the last using
a small change. The standard minimal-change algorithms for permutation
generation are given  in \cite[7.2.1.2]{knuth;taocp4}. A very general method of 
enumerating permutations is as follows. A \Em{Sims table} for a 
subgroup $G$ of $\sym_n$ is a family of 
subsets $S_1, \dots $ of $G$ having the following property: for each $j, k$ with 
$1\leq j \leq k \leq n$, $S_k$ contains exactly one element that fixes all elements greater than 
$k$ and takes $k$ to $j$, whenever $G$ itself contains such a permutation. It is easily seen 
\cite[Lemma S]{knuth;taocp4} that if $S_1, \dots, S_{n-1}$ is a Sims table then every element 
of $G$ has a unique representation as a product $\perm = \perm_1 \dots \perm_{n-1}$, where 
$\perm_k \in S_k$ for each $k$. There is also a  unique dual representation of the form 
$\perm = \perm_{n-1}^{-1} \dots \perm_1^{-1}$ with $\perm_k \in S_k$, obtained by inverting the first
representation for $\perm^{-1}$. 

An inspection of the proofs shows that $G$ need not be a group for such results to hold. In fact, 
it is only necessary that $G$ be closed under taking inverses. Thus, for example, the set $\cyc_n$ 
could be used. 

The triangular decomposition fits into this framework. For each $k$, the set $S_k$ consists of 
all transpositions $\trans_{jk}$ with $j\leq k$ (in the case $\sym_n$) or $j < k$ 
(in the case $\cyc_n$). The Sims representation with respect to these sets $S_k$ 
is precisely the triangular representation. 

\begin{table}
\begin{tabular}{|cccc|cccc|}
\hline
\multicolumn{4}{|c|}{Lex order on $[n] \times \dots \times [1]$} & 
\multicolumn{4}{|c|}{Induced order on $\sym_n$} \\
\hline
0000 & 0100 & 0200 & 0300 & 1230 & 3201 & 1302 & 1203 \\  
0001 & 0101 & 0201 & 0301 & 2130 & 2301 & 3102 & 2103 \\
0010 & 0110 & 0210 & 0310 & 2310 & 2031 & 3012 & 2013 \\
0011 & 0111 & 0211 & 0311 & 3210 & 0231 & 0312 & 0213 \\
0020 & 0120 & 0220 & 0320 & 1320 & 3021 & 1032 & 1023 \\
0021 & 0121 & 0221 & 0321 & 3120 & 0321 & 0132 & 0123 \\
\hline
\end{tabular}
\caption{The Fisher-Yates encoding for $n=4$.}
\label{t:FY encode}
\end{table}

\begin{table}
\begin{tabular}{|cccccc|cccccc|}
\hline
\multicolumn{6}{|c|}{Lex order on $[1] \times \dots \times [n]$} & 
\multicolumn{6}{|c|}{Induced order on $\sym_n$} \\
\hline
0000 & 0010 & 0020 & 0100 & 0110 & 0120 & 1230 & 2310 & 1320 & 2130 & 3210 & 3120 \\  
0001 & 0011 & 0021 & 0101 & 0111 & 0121 & 3201 & 2031 & 3021 & 2301 & 0231 & 0321 \\
0002 & 0012 & 0022 & 0102 & 0112 & 0122 & 1302 & 3012 & 1032 & 3102 & 0312 & 0132 \\
0003 & 0013 & 0023 & 0103 & 0113 & 0123 & 1203 & 2013 & 1023 & 2103 & 0213 & 0123 \\
\hline
\end{tabular}
\caption{The Fisher-Yates encoding for $n=4$.}
\label{t:FY dual encode}
\end{table}

For comparison we include the usual inversion encoding in Table~\ref{t:inv encode}.

\begin{table}
\begin{tabular}{|cccccc|cccccc|}
\hline
\multicolumn{6}{|c|}{Lex order on inversion function} & 
\multicolumn{6}{|c|}{Induced order on $\sym_n$} \\
\hline
0000 & 0010 & 0020 & 0100 & 0110 & 0120 & 0123 & 0213 & 2013 & 1023 & 1203 & 2103 \\  
0001 & 0011 & 0021 & 0101 & 0111 & 0121 & 0132 & 0231 & 2031 & 1032 & 1230 & 2130 \\
0002 & 0012 & 0022 & 0102 & 0112 & 0122 & 0312 & 0321 & 2301 & 1302 & 1320 & 2310 \\
0003 & 0013 & 0023 & 0103 & 0113 & 0123 & 3012 & 3021 & 3201 & 3102 & 3120 & 3210 \\
\hline
\end{tabular}
\caption{The inversion encoding for $n=4$.}
\label{t:inv encode}
\end{table}

A \Em{Gray code} for $G = \sym_n$ is a Hamiltonian path in the Cayley graph of $G$ where the 
generating set is the set of all transpositions. The usual 
Gray code on words in $[1] \times \dots \times [n]$
induces a Gray code on $G$ via the inversion encoding, since each minimal 
change to a word corresponds to a transposition of adjacent symbols. This can be seen by reading 
the columns alternately downwards and upwards from left to right in Table~\ref{t:inv encode}.

What happens when we use instead the Fisher-Yates encoding? The Gray code order on words 
induces a Hamiltonian path in the Cayley graph of $G$, but with respect to a different set 
of generators. The generators
in question are in fact transpositions and $3$-cycles. To see this, note that to get from one entry 
to the next we move from $\perm_1 \trans \perm_2$ to $\perm_1 \trans' \perm_2$, which is achieved 
by multiplying by $\perm_2^{-1} \trans' \trans \perm_2$. Since $\trans$ and $\trans'$ transpose 
some symbol 
$k$ with $j, j'$ respectively, where $j \neq j'$ and $j, j' \leq k$, the product $\trans' \trans$ is 
the permutation that moves $j'$ to $j$, $j$ to $k$ and $k$ to $j'$. If either $j = k$ or $j' = k$ then
$\trans' \trans$ is a transposition, and otherwise it is a $3$-cycle. The conjugation by $\perm_2$ 
preserves the cycle structure.
 
Note that the restriction to $\cyc$ is better behaved and the Cayley graph of the set $\cyc$ 
with respect to the set of 3-cycles has a Hamiltonian cycle. Since we always have $j, j' < k$ in the 
case of $\cyc$, the transpositions are never needed and we always move from one element to the next by 
multiplying by a $3$-cycle. Since smaller changes could only be transpositions, and multiplying 
an $n$-cycle by a transposition can never yield an $n$-cycle, the enumeration described above 
deserves the name ``\Em{Gray code for cycles}". 

\section{Analysis of some quantities}
\label{sec:analysis}

Obvious quantities to be studied are: the number of swaps; the number of
times a given symbol is chosen by the random calls (we call this the number of moves, although 
some of these moves will be trivial); the total distance moved by a given symbol; 
the total distance moved. The second and third of these were discussed in
\cite{prodinger;sattolo, mahmoud;sattolo, wilson;sattolo} for the case of Sattolo's algorithm.

The number of swaps is always $n - 1$ for each algorithm, but some of these
can be trivial (and hence executed more quickly) for the Fisher-Yates algorithm, whereas
every exchange is nontrivial in Sattolo's algorithm. The number of
nontrival swaps is the number of elements moved by the permutation, or $n -
f$, where $f$ is the number of fixed points. The generating function for
permutations by size and fixed points is well known to be 
$$
\sum_{\perm} \frac{x^{|\perm|}}{|\pi|!} u^{f(\pi)} = 
\frac{e^{(u - 1) z}}{1 - z}.
$$ 
For example, the expected number of fixed points is $1$ for every $n$.

\subsection*{Number of moves and distance moved by an element}
To avoid excessive case distinctions we consider the slight variant of these algorithms in which the 
final ``swap" of $\perm[0]$ with itself is performed (this corresponds to the ``downto" loops 
in Figure~\ref{fig:pseudo} going down to $0$ instead of $1$. 

We consider normalized counting generating functions of the form 

\begin{equation*}
F(u, t, x) := \sum_{\perm \in \sym, p \in [n(\perm)]} u^{\chi(\perm, p)}
t^p \frac{x^{n(\perm)}}{|\sym_{n(\perm)}|} 
= \sum_{n\geq 1} \frac{x^n}{n!} \sum_{1\leq p \leq n} t^p 
\sum_{\perm\in\sym_n} u^{\chi(\perm, p)}.
\end{equation*}

An auxiliary ``diagonal" GF will also be useful:

\begin{equation*}
G(u, x) := \sum_{\perm\in \sym} u^{\chi(\perm, n(\perm))} 
\frac{x^{n(\perm)}}{|\sym_{n(\perm)}|} 
= \sum_{n\geq 1} \frac{x^n}{n!} \sum_{\perm\in\sym_n} u^{\chi(\perm, n)}.
\end{equation*}

Here $\chi$ is a given parameter of interest such as number of moves, etc.
Of course $F$ and $G$ can be interpreted probabilistically as ``grand"
PGFs. For example, if $\chi(\perm, p)$ is the number of moves made by $p$
in obtaining $\perm$ via the Fisher-Yates algorithm, and $M_{np}$ the
random variable obtained by evaluating $\chi$ at an element of $\cyc_n$
chosen uniformly at random, then letting $\phi_{np}(u) = \sum_{l\geq 0}
\mathbb{P}(M_{np} = l) u^l$ denote the PGF of  $M_{np}$, we have

$$
F(u, t, x) = \sum_{n\geq 1} x^n \sum_{p=1}^n t^p \phi_{np}(u).
$$

We first consider the case where $\chi$ is the number of moves of a given symbol. 
The triangular decomposition yields the recurrence

\begin{equation}
\label{eq:chirec1}
\chi(\perm, p) = 
\begin{cases}
\chi(\down{\perm}, p) & \text{if $p\neq n(\perm), p\neq q(\perm)$;} \\
1 + \chi(\down{\perm}, q(\perm))  &\text{if $p = n(\perm), 
p\neq q(\perm)$;} \\
1 & \text{if $p\neq n(\perm), p = q(\perm)$;} \\
1 & \text{if $p = n(\perm), p = q(\perm)$.} \\
\end{cases}
\end{equation}

In the case where $\chi$ is the distance moved by an element, we have the recurrence
\begin{equation}
\label{eq:chirec2}
\chi(\perm, p) = 
\begin{cases}
\chi(\down{\perm}, p) & \text{if $p \neq  n(\perm), p\neq q(\perm)$;} \\
n(\perm) - q(\perm) +  \chi(\down{\perm}, q(\perm)) 
&\text{if $p = n(\perm), p \neq q(\perm)$;} \\
n(\perm) - q(\perm) & \text{if $p \neq n(\perm), p =  q(\perm)$;} \\
0 & \text{if $p = n(\perm), p = q(\perm)$.} \\
\end{cases}
\end{equation}

We partition the index set $\indset = \{ (\perm, p) \mid \perm \in
\sym, 1\leq p \leq n(\perm) \}$ into $4$ disjoint subsets $\indset_1,
\dots, \indset_4$ according to the cases just listed. Denote by
$\Sigma_k(u, t, x)$ the part of the sum  defining $F$ corresponding to
index set $\indset_k$, so that $F = \Sigma_1 + \Sigma_2 + \Sigma_3 + \Sigma_4$.

Note that for each $\chi$ we have
\begin{align*}
\Sigma_4 + \Sigma_2 & = \sum_{n\geq 1} \frac{x^n}{n!} t^n 
\sum_{\perm\in\sym_n} u^{\chi(\perm, n )} \\
& = G(u, tx).
\end{align*}

In the sum $\Sigma_1$, indices $p\in [n(\perm)]$ satisfying the conditions
$p\neq n(\perm), p \neq q(\perm)$ occur if and only if $ n(\perm) \geq 2$. The  set 
$\indset_1$ is in bijection with the set
$$
\{(\perm, q, p ) \mid n(\perm) \geq 2, 1\leq p < n(\perm), 1 \leq q
\leq n(\perm), p \neq q\}.
$$ 

Let $A(u, t, x)$ be the antiderivative of $F(u, t, x)$ with respect to $x$ having 
$A(u,t,0) = 0$. Then for each $\chi$ we obtain

\begin{align*}
\Sigma_1(u, t, x) & = 
\sum_{n\geq 2}\frac{x^{n}}{n!} \sum_{\perm\in \sym_n} 
\sum_{1\leq p < n, p\neq q(\perm)}  t^p   u^{\chi(\perm, p)}  \\
& = \sum_{n\geq 2} \frac{x^n}{n!} \sum_{\perm\in \sym_{n-1}}
\sum_{1\leq q \leq n} \sum_{1 \leq p < n, p \neq q} 
u^{\chi(\perm, p)} t^p \\
& = x \sum_{n\geq 1} \frac{x^n}{(n+1)!} \sum_{1 \leq p \leq n} t^p 
\sum_{\perm\in \sym_n} u^{\chi(\perm, p)} \sum_{1\leq q \leq n+1, q \neq p} 1\\
& = x \sum_{n\geq 1} \frac{n x^n}{(n+1)!} \sum_{1\leq p \leq n} 
t^p \sum_{\perm\in \sym_n}  u^{\chi(\perm, p)}  \\
& = x \sum_{n\geq 1} \frac{x^n}{n!} \sum_{1\leq p \leq n} 
t^p \sum_{\perm\in \sym_n}  u^{\chi(\perm, p)} - \sum_{n\geq 1} \frac{ x^{n+1}}{(n+1)!} \sum_{1\leq p \leq n} 
t^p \sum_{\perm\in \sym_n}  u^{\chi(\perm, p)} \\
& = x F(u, t, x) - A(u, t, x).
\end{align*}

We now determine $\Sigma_3(u, t, x)$. The set $\indset_3$ is in bijection with 
$\{\perm \in \sym \mid n(\perm) \geq 2, q \neq n(\perm)\}$. Thus for the number of moves we have
\begin{align*}
\Sigma_3(u, t, x) 
& = \sum_{n\geq 2} \frac{x^n}{n!} \sum_{\perm\in\sym_n} t^{q(\perm)} u^1 
 = u \sum_{n\geq 2} \frac{x^n}{n!} \sum_{\down{\perm}\in\sym_{n-1}} 
\sum_{q=1}^{n} t^q\\
& = u \sum_{n\geq 1} \frac{x^{n+1}}{(n+1)!} \sum_{\perm\in\sym_{n}} 
\sum_{q=1}^{n} t^q
 = u\sum_{n\geq 1} \frac{x^{n+1}}{n+1} \sum_{q=1}^{n} t^q\\
& = \frac{utx}{1-t} \sum_{n\geq 1} \frac{x^{n} (1 - t^{n})}{n+1}
 = \frac{u}{1-t} [\log(1 - tx)- t \log(1 - x)].
\end{align*}

Note that when $t=1$ we have the formula $\Sigma_3(u,1,x) = u \log(1-x) + ux/(1-x)$.

To obtain $\Sigma_3$ for the distance moved, a similar calculation yields 
$$
\frac{u}{u - t} \left[u \log(1 - tx) - t \log(1 - ux)\right].
$$
 
We now consider $G$. We have for the number of moves
\begin{align*}
G(u, x) & = \sum_{n\geq 1} \frac{x^n}{n!} \sum_{\perm\in\sym_n} 
u^{\chi(\perm, n)} 
 = \sum_{n\geq 2} \frac{x^n}{n!} \sum_{\perm\in\sym_n, q(\perm) \neq n(\perm)} 
u^{\chi(\perm, n)} + \sum_{n\geq 1} \frac{x^n}{n!} \sum_{\perm\in\sym_n, q(\perm) = n(\perm)} 
u^{\chi(\perm, n)} \\
& = \sum_{(\down{\perm}, 1\leq q \leq n(\down{\perm}))} u^{1 + \chi(\down{\perm}, q)} 
\frac{x^{n(\down{\perm}) + 1}}{n(\down{\perm} + 1)!} + \sum_{n\geq 1}\frac{x^n}{n!} 
\sum_{\perm\in\sym_n, q(\perm) = n(\perm)} u^1  \\
 & =  u \sum_{(\perm, 1\leq q \leq n(\perm))} u^{\chi(\perm, q)} \frac{x^{n(\perm) + 1}}{(n(\perm) + 1)!} 
 + u \sum_{n\geq 1} \frac{x^n}{n!} (n-1)! \\
& =  u A(u, 1, x) - u \log(1 - x).
\end{align*}
Similarly in the case of distance moved we obtain $G(u, x) = A(u, u^{-1}, ux) - \log(1 - x)$.

Thus for the number of moves, by differentiating we obtain the system

\begin{align*}
(1 - x) F'(u, t, x) & = t G'(u, tx) + \Sigma'_3(u, t, x) \\
G'(u, x) & =  u F(u, 1, x) + \frac{u}{1-x}.
\end{align*}

Substituting $t = 1$ and eliminating $G'(u,x)$ we obtain 
\begin{align*}
(1 - x) F'(u, 1, x) - uF(u, 1, x) & = \frac{u}{1 - x} + \Sigma'_3(u, 1, x) 
 = \frac{u}{(1 - x)^2} \\
F(u, 1, 0) & = 0
\end{align*}

which  yields

\begin{align*}
F(u, 1, x) & = \frac{u}{2 - u} \left[(1 - x)^{-2} - (1 - x)^{-u}\right].
\end{align*}
From this $G$ can be found explicitly via a single integration.

$$
G(u, x) = \frac{u^2}{2 - u} \left[ \frac{1}{1 - x} + \frac{(1 - x)^{1 - u}}{1 - u}\right] 
- \frac{u^2}{1-u} - u \log (1 - x).
$$

To find $F$ explicitly is more difficult, because it requires the integration of $(1-tx)^{-u}(1-x)^{-1}$, 
and we do not pursue it here. In any case we have the defining equation
\begin{align*}
(1 - x)F'(u, t, x) & = utF(u, 1, tx) + \frac{u}{1 - tx} + \Sigma'_3(u, t, x) \\
F(u, t, 0) & = 0.
\end{align*}

Using this defining equation we may easily extract the coefficient of 
$x^n t^p$ to obtain the probability generating function $\phi_{np}(u)$, 
or extract moments by evaluating appropriate partial $u$-derivatives at $u = 1$ as usual. 
For example, $\phi_{np}(u)$ is the coefficient of $x^n t^p$ in $F$, hence
equals $(1/n)$ times the coefficient of $x^n t^p$ in $xF'$, and the mean of
the random variable with PGF $\phi_{np}(u)$ can therefore be obtained by
evaluating $x \partial^2 F/\partial x \partial u$ at $u = 1$, then dividing by
$n$. The defining equation allows us to express these derivatives in terms of derivatives of the known series 
$G(u, x)$ and  $F(u, 1, x)$. 

In detail, we see that (with subscripts denoting partial derivatives, and 
$\Sigma = \Sigma_3$ to avoid notational overload)
\begin{align*}
x(1-x) F_{13}(1, t, x) & = txG_{12}(1, tx) + x\Sigma_{13}(1, t, x) \\
& = \frac{tx}{(1-tx)^2} + \frac{2t^2x^2}{(1-tx)^2} + \frac{tx \log (1 - tx)}{1 - tx} + 
\frac{tx}{(1-x)(1-tx)}.
\end{align*}
Extracting the coefficient of $t^p$ from the right side yields $px^p + 2(p-1) x^p  + x^{p+1}/(1-x) - H_{p-1}$ 
where $H_{p}$ denotes the $p$th harmonic number $\sum_{1\leq i\leq p} 1/i$.
Dividing by $(1-x)$ and extracting the coefficient of $x^n$ yields
(where $M$ denotes the number of moves)
$$
E[M_{np}]  = \frac{n + 2p - 2 - H_{p-1}}{n}.
$$

Higher moments can also be obtained with more calculation of the same type, but we do not pursue this 
aspect here.

We can also immediately extract recurrences for the probability generating function. We obtain 
\begin{align*}
\phi_{np}(u) & = \frac{p}{n}\phi_{pp}(u) + (1 - p/n) u \\
\phi_{nn}(u) & = \frac{u}{n} \left[ 1 + \sum_{1\leq p \leq n-1} \phi_{n-1, p}(u) \right]. 
\end{align*}
We can now easily write down an explicit formula for the probability generating 
function $\phi_{np}$ by extracting of coefficients from $G$ and using the first recurrence above.
We have
\begin{align*}
\phi_{nn}(u) & = \frac{u}{n} + \frac{u^2}{2-u} \left[1 - \frac{u(u+1)\dots (u+n-2)}{n!}\right] \\
\phi_{np}(u) & = \frac{n-(p-1)}{n} u + \frac{p}{n} \frac{u^2}{2-u} \left[1 - \frac{u(u+1)\dots (u+p-2)}{p!} \right].
\end{align*}

Similarly for the distance moved we obtain

\begin{align*}
(1 - x) F'(u, t, x) & = t G'(u, tx) + \Sigma'_3(u, t, x) \\
G'(u, x) & =  uF(u, u^{-1}, ux) + (1-x)^{-1}
\end{align*}

which leads via the substitution $t \leftarrow u^{-1}$ to 
\begin{align*}
(1 - x) F'(u, u^{-1}, x) - F(u, u^{-1}, x)  & =  (1 - u^{-1}x)^{-1} + \Sigma'_3(u, u^{-1}, x) \\
F(u, u^{-1}, 0) & = 0.
\end{align*}

This equation is exact and leads to 
\begin{align*}
F(u, u^{-1}, ux) & =  \frac{1}{1-ux} \left[- u \log (1 - x) + \Sigma_3(u, u^{-1}, ux) \right] \\
& =  \frac{1}{1-ux} \left[ - u \log (1 - x)
+ \frac{u}{1 - u^2} \left[ \log(1 - u^2x) - u^2\log(1 - x) \right] \right].
\end{align*}

Another integration yields $G(u, x)$ and again we have a defining equation for $F(u, t, x)$. Each of these 
integrations takes us outside the realm of elementary functions. However, we can 
systematically extract coefficients as before from the equations that we have. We omit any details of the 
calculations. The probability generating function $\xi_{np}$ of the random variable 
$D_{np}$ satisfies the recurrence
\begin{align*}
\xi_{np}(u) & = \frac{p}{n}\xi_{pp}(u) + \frac{1}{n} \sum_{j=1}^{n-p} u^{j}\\
\xi_{nn}(u) & = \frac{1}{n} \left[ 1 + \sum_{p<n} u^{n-p} \xi_{n-1,p}(u) \right].
\end{align*}
They have explicit formulae such as
\begin{align*}
\xi_{nn}(u) & = \frac{1}{n} \left[1 + \sum_{j=1}^{n-1} \frac{u^{n-j}}{j} +
\frac{u^{n+2}}{1 - u^2} \sum_{j=1}^{n-1} \frac{u^{1 - j} - u^{j - 1}}{j} \right]. 
\end{align*}
Similarly we could extract the mean and higher moments as before. From the explicit forms one could consider 
limit distributions by considering the pointwise limit of the probability generating function $\phi_{np}$ 
or $\xi_{np}$. We do not pursue this further here as the computations are routine but tedious.

\subsection*{The total distance moved}
\label{ss:dist}

We consider the total distance moved rightwards by elements, $D_n =
\sum_p D_{np}$ (this of course equals the total distance moved leftwards
by elements). As a random variable, $D_n$ is the sum of $D_{n-1}$ and a
random variable $U_i$ that is uniform on $[0..n-1]$. Thus $D_n$ is
distributed as $\sum_{i=2}^n U_i$ and has probability generating
function $\prod_{i=1}^n \frac{1 - u^i}{1 - u}$.

Note that this PGF is the same as the one for inversions. Thus the
number of inversions and the total rightward distance have the same
distribution (in other words, $D_n$ is a \Em{Mahonian statistic}). Hence for
each $k$, the number of permutations in $\sym_n$ with $k$ inversions is the
same as the number of permutations in $\sym_n$ whose rightward distance is
$k$. 

\section{Extensions and discussion}
\label{sec:extensions}

Despite an extensive literature search, I can only find two places in the literature in which the 
very natural Fisher-Yates encoding is mentioned. In neither paper was the connection with 
cyclic permutations mentioned. 

What we have called the Fisher-Yates encoding was used in \cite{MaRa2001} to study anti-excedances.
In \cite{myrvold-ruskey}, two unranking and ranking functions for
permutations were presented, each taking linear time to compute. Although
not mentioned in that paper, it is easily seen that those orderings
correspond via the Fisher-Yates encoding to lexicographic order on the triangular cartesian product 
$[1]\times \dots \times [n]$ or $[n] \times \dots \times [1]$ as we have described above.

We note that Sattolo's algorithm is a special case of a method to uniform generation of 
permutations with a fixed number $k$ of cycles \cite{wilf-eastwest}.

Suppose that $m \geq 0$ is fixed, and at the $i$th step, we swap $\perm(i)$
and $\perm(j)$ where $j$ is chosen uniformly from $[1..n+1-i-m]$. The cases 
$m = 0$ and $m = 1$ respectively correspond to the Fisher-Yates and Sattolo
algorithms. Other values of $m$ do not appear to be particularly interesting, although 
we have not pursued them.

I have not yet been able to determine whether the statistic $D$ above is known.
An interesting question is its correlation with other well-known permutation statistics. 

There is a small connection between the Fisher-Yates algorithm and sorting. 

\begin{proposition}
Let $\perm \in \sym_n$ and let $\trans_1 \dots \trans_{n-1}$ be its 
triangular decomposition. Then in order to sort $\perm^{-1}$, selection
sort applies the transpositions $\trans_{n-1}, \dots , \trans_{1}$ in that
order. 
\end{proposition} 
\begin{proof}
Selection sort first chooses $n$ and puts it in the correct position; this corresponds to
 postmultiplication by $\trans(n, n  \perm)$, which corresponds to premultiplication by 
 $\trans(n, n  \perm^{-1})$. The result follows by induction.
\end{proof}

For example, to generate $\perm = 3 4 2 1$ the algorithm proceeds as
follows: $1 2 3 4$, $4 2 3 1$, $4 3 2 1$, $3 4 2 1$, yielding the strict
triangular decomposition $(1 2)(2 3)(1 4)$. Applying these in turn to
$\perm^{-1}$ gives $4 3 1 2, 2 3 1 4, 2 1 3 4, 1 2 3 4$ which is the
list created by selection sort when sorting $\perm^{-1}=4 3 1 2$.


\begin{thebibliography}{Mah03}
\providecommand{\bysame}{\leavevmode\hbox to3em{\hrulefill}\thinspace}

\bibitem[Durs1964]{durstenfeld} Richard Durstenfeld, 
\emph{Algorithm 235: Random permutation}, Comm. Assoc. Comput. Mach. 7, 1964, 420.

\bibitem[FY1938]{fisher-yates} R. A. Fisher and F. Yates, Example 12, 
\emph{Statistical Tables}, London, 1938.

\bibitem[GX1988]{gries-xue}
David Gries and Jin~Yun Xue, 
\emph{Generating a random cyclic permutation}, BIT \textbf{28} (1988),  569--572. 
  
\bibitem[Knut1969]{knuth;taocp2} Donald E. Knuth, 
\emph{The art of computer programming. Vol. 2: Seminumerical algorithms.} 
Addison-Wesley, 1969.
  
\bibitem[Knut2004]{knuth;taocp4}
Donald E. Knuth, \emph{The art of computer programming. Vol. 4, Fasc. 2.
Generating all tuples and permutations.} Addison-Wesley, 2005. 

\bibitem[MaRa2001]{MaRa2001}
Roberto Mantaci and Fanja Rakotondrajao, \emph{A permutation representation that knows what 
``Eulerian" means}, Discrete Math. Theor. Comput. Sci. \textbf{4} (2001), 101--108.

\bibitem[Mahm2003]{mahmoud;sattolo}
Hosam~M. Mahmoud, 
\emph{Mixed distributions in {S}attolo's algorithm for cyclic permutations via 
randomization and derandomization}, 
J. Appl. Probab. \textbf{40} (2003), 790--796. 
 
\bibitem[MyRu2001]{myrvold-ruskey} Wendy Myrvold and Frank Ruskey, 
\emph{Ranking and unranking permutations in linear time}, 
Inform. Process. Lett. 79 (2001), 281--284. 

\bibitem[Prod2002]{prodinger;sattolo}
Helmut Prodinger, 
\emph{On the analysis of an algorithm to generate a random cyclic permutation},
 Ars Combin. \textbf{65} (2002), 75--78. 

\bibitem[Prod]{prodinger;webnote} 
\bysame, Online document at \newline
\texttt{http://math.sun.ac.za/\~{}prodinger/abstract/abs\_161.htm}.

\bibitem[Satt1986]{sattolo}
Sandra Sattolo, 
\emph{An algorithm to generate a random cyclic permutation}, 
Inform. Process. Lett.\textbf{22} (1986), 315--317. 

\bibitem[Wilf]{wilf-eastwest}
Herbert Wilf,
\emph{East Side, West Side}, 
lecture notes available from \texttt{http://www.cis.upenn.edu/\~{}wilf/lecnotes.html}.

\bibitem[Wils2004]{wilson;sattolo}
Mark C. Wilson, 
\emph{Probability generating functions for Sattolo's algorithm}, 
J. Iranian Stat. Soc. 3 (2004), 297--308.

\end{thebibliography}
\end{document}